\documentclass[a4]{amsart}
\usepackage{graphicx} 
\usepackage{amssymb,amsmath, amsthm}
\usepackage{enumerate}
\usepackage{cite}
\usepackage{soul}
\usepackage{comment}
\usepackage{amsfonts}
\usepackage{xcolor}
\usepackage{url}
\usepackage{hyperref}

\usepackage{doi}        
\usepackage[numbers]{natbib}  

\newtheorem{Thm}{Theorem}[section]
\newtheorem{Lem}[Thm]{Lemma}
\newtheorem{rem}[Thm]{Remark}
\newtheorem{prop}[Thm]{Proposition}
\newtheorem{Cor}[Thm]{Corollary}

\newcommand{\Four}{\mathcal{F}}

\newcommand{\N}{\mathbb{N}}
\newcommand{\Sc}{\mathcal{S}}

\newcommand{\eps}{\varepsilon}
\newcommand{\Class}{\mathcal{C}}
\newcommand{\Haus}{\mathcal{H}} 
\newcommand{\Prob}{\mathcal{P}}

\newcommand{\Rsp}{\mathbb{R}}
\newcommand{\Sph}{\mathcal{S}} 

\newcommand{\abs}[1]{\left\vert #1 \right\vert}
\newcommand{\sca}[2]{\langle #1 |#2\rangle}
\newcommand{\dd}{\mathrm{d}}
\newcommand{\Lip}{\mathrm{Lip}}
\newcommand{\Expect}{\mathbb{E}}
\newcommand{\Wass}{\mathcal{W}}
\newcommand{\SlicedWass}{\mathcal{SW}}
\newcommand{\MaxSlicedWass}{\mathcal{MSW}}
\newcommand{\hide}[1]{}

 \newcommand{\edz}[1]{\sideremark{#1}}

\def\sideremark#1{\ifvmode\leavevmode\fi\vadjust{\vbox to0pt{\vss
 \hbox to 0pt{\hskip\hsize\hskip1em
 \vbox{\hsize3cm\tiny\raggedright\pretolerance10000
 \noindent #1\hfill}\hss}\vbox to8pt{\vfil}\vss}}}

\title[Sliced and standard Wasserstein distances]{Sharp comparisons between sliced and standard $1$-Wasserstein distances}

\author{Guillaume Carlier}
\address{CEREMADE, UMR CNRS 7534, Universit\'e Paris
Dauphine, PSL, Pl. de Lattre de Tassigny, 75775 Paris Cedex 16, France \and INRIA-Paris}
\email{carlier@ceremade.dauphine.fr}

\author{Alessio Figalli}
\address{ETH Z\"{u}rich, Department of Mathematics, R\"{a}mistrasse 101, 8092, Z\"{u}rich, Switzerland}
\email{alessio.figalli@math.ethz.ch}

\author{Quentin M\'erigot}
\address{Universit\'e Paris-Saclay, CNRS, Laboratoire de math\'ematiques d'Orsay, 91405, Orsay, France \and Institut universitaire de France}
\email{quentin.merigot@universite-paris-saclay.fr}

\author{Yi Wang}
\address{Johns Hopkins University, 404 Krieger Hall, Department of Mathematics, 3400 N. Charles St, Baltimore, MD 21218}
\email{ywang261@jhu.edu}

\date{\today}

\begin{document}

\maketitle

\begin{abstract} Sliced Wasserstein distances are widely used in practice  as a computationally efficient alternative to Wasserstein distances in high dimensions. In this paper, motivated by theoretical foundations of this alternative, we prove quantitative estimates between the sliced $1$-Wasserstein distance and the $1$-Wasserstein distance. We construct a concrete example to demonstrate the exponents in the estimate is sharp. We also provide a general analysis for the case where slicing involves projections onto $k$-planes and not just lines.
\end{abstract}

\section{Introduction}

 The sliced Wasserstein distance of order  $p\geq 1$ between two probability measures $\mu$ and $\nu$ on $\mathbb{R}^d$
with finite $p$-moments, $\SlicedWass_p(\mu,\nu)$ is obtained by  averaging the  (one-dimensional) 
Wasserstein distances between their orthogonal projections onto lines:
\begin{equation}\label{eq:def-sliced}
\SlicedWass_p(\mu,\nu) :=\left( \int_{\theta \in\Sph^{d-1}}
\Wass_p(P_{\theta\#}\mu,P_{\theta\#}\nu)^p \dd \sigma_{d-1} ( \theta)\right)^{\frac{1}{p}}
\end{equation}
where $\Sph^{d-1}$ is the unit sphere of $\Rsp^d$, $\sigma_{d-1}=k_d \Haus^{d-1}\llcorner \Sph^{d-1}$ is the uniform probability measure on $\Sph^{d-1}$ and $P_\theta(x) = \sca{\theta}{x}$. 
Sliced Wasserstein distances have originally been introduced by Marc
Bernot around 2007 as a computationally inexpensive variant of  optimal
transport distances. Indeed, the first obvious appeal of sliced Wasserstein distances is their ability to replace the resolution of (difficult) optimal transport problems in dimension $d\geq 2$ with a family of unidimensional problems which can be solved in (almost) closed form using standard monotone transport constructions. In addition, these unidimensional transport distances can be computed independently, enabling efficient parallelization.  These properties have made sliced Wasserstein distances increasingly popular in a range of applied domains.
They were first used in texture synthesis and, more generally, in imaging science~\cite{rabin2012wasserstein}, and have since become a valuable alternative to Wasserstein distances in many applications in high-dimensional statistics and machine learning~\cite{Kolouri2019a, Kolouri2019b, Meunier2022, Rustamov2020, Bonet2023}. We note, however, that the aforementioned computational advantage is partly offset by the need to compute the high-dimensional integral in the definition of $\SlicedWass_p$, which remains a nontrivial numerical challenge and has been the focus of recent work \cite{Bonneel2015}, \cite{nadjahi2021sliced}, \cite{nadjahi2020statistical}. 
Another popular approach to alleviate the computational burden of Wasserstein distances is entropic regularization and the associated Sinkhorn algorithm, popularized by Cuturi~\cite{Cuturi2013}. 

\smallskip

Another advantage of sliced Wasserstein distances is that they exhibit a nearly dimension-independent statistical behaviour. For instance,  \cite[Corollary 2]{nadjahi2020statistical}, based on \cite{fournier2015rate}, implies that for compactly supported probability measures $\mu,\nu$ on $\Rsp^d$, one has 
$$ \Expect\abs{\SlicedWass_p(\mu_N,\nu_N) - \SlicedWass_p(\mu,\nu)} \lesssim \left(\frac{1}{N}\right)^{\frac{1}{2p}},$$
where $\mu_N,\nu_N$ denote the empirical measures associated to i.i.d. $N$-samples drawn from $\mu$ and $\nu$ respectively.

\smallskip

Many variants of sliced Wasserstein distances have been considered, among them the max-sliced distances \cite{deshpande2019max}, where the integral in equation~\eqref{eq:def-sliced} is replaced by a maximum
\begin{equation}\label{eq:def-max-sliced}
\MaxSlicedWass_p(\mu,\nu) := \max_{\theta\in\Sph^{d-1}} \Wass_p(P_{\theta\#}\mu, P_{\theta\#}\nu).
\end{equation}
This is a special case  ($q=\infty$) of the family  of sliced distances  studied in \cite{kitagawa2024}, which involve two parameters: $q$ for the $L^q$ norm with respect to the $\theta$ variable and $p$ for the exponent of the Wasserstein distance. It is worth noting that the comparison between  $\MaxSlicedWass_1$ and $\Wass_1$ was already addressed  40 years ago by Hahn and Quinto in the interesting (though likely not widely known within the optimal transport community) article \cite{Hahn-Quinto}. This result was very recently improved by Bobkov and G\"{o}tze in \cite{bobkov2024}. The fine study of the geometry and analytical and topological properties of spaces of probability measures over $\Rsp^d$ equipped with sliced Wasserstein distances, as  metric spaces, is the object of a very active stream of recent research, \cite{park2023geometry}, \cite{kitagawa2024}. 

\smallskip

As far as the comparison between the sliced distance $\SlicedWass_p$ and  the usual Wasserstein distance $\Wass_p$ is concerned, the first obvious observation is that $\SlicedWass_p \leq \Wass_p$. This follows from the $1$-Lipschitz property of projections. In addition, if we restrict ourselves to probabilities supported on a fixed closed ball, it is easy to see that these two distances metrize weak$^*$ convergence, hence are \emph{topologically} equivalent (see \cite{kitagawa2024} for the extension of this topological equivalence to the whole of $\Rsp^d$). On the other hand, it is well-known that $\SlicedWass_p$ and $\Wass_p$ are not \emph{metrically} equivalent.

To see this, fix $d=2$  and let $X$ and $Y$ be two independent standard Gaussian random variables on $\mathbb R$. 
Fix a small $\eps>0$, and define $\mu$ as the law of $(X,0)$ and $\nu_\eps$ as the law of $(X,\eps Y)$.
Then $\Wass_1(\mu,\nu_\eps)=\eps \Expect(|Y|)$. For every $\theta=(\cos(\alpha), \sin(\alpha))\in \Sph^{1}$, the projections ${P_\theta}_\# \mu$  and ${P_\theta}_\# \nu_\eps$ are Gaussians with respective standard deviations $\vert \cos (\alpha)\vert$ and $(\cos^2(\alpha)+\eps^2 \sin^2(\alpha))^{1/2}$.
Therefore, $\SlicedWass_1(\mu, \nu_\eps)$ has an explicit expression:
\begin{multline*}\SlicedWass_1(\mu, \nu_\eps)= \frac{\Expect(\vert Y\vert) }{2\pi} \int_0^{2 \pi}  \Big(\sqrt{\cos^2(\alpha)+\eps^2 \sin^2(\alpha)}-\vert \cos(\alpha)\vert \Big) \dd \alpha \\
\lesssim \eps^2 \vert \log \eps | \simeq \Wass_1(\mu, \nu_\eps)^2\big|\log \Wass_1(\mu, \nu_\eps)\big|.
\end{multline*}
This example shows that there is no bound of the form $\Wass_1\leq C \SlicedWass_1^\alpha$ for $\alpha>\frac{1}{2}$.
\smallskip

Quantitative control of $\Wass_p$ by  $\SlicedWass_p$ (or $\MaxSlicedWass_p$) and how it is affected by the curse of dimensionality  is a key issue for justifying the practical use of sliced distances as a reliable (or not) surrogate of $\Wass_p$ in high dimensions.  In this  direction, let us quote two  results for $p=1$ and $\mu$ and $\nu$ supported on the $d$-dimensional ball $B_R$ of radius $R$.  The first one has been established in the PhD thesis of Nicolas Bonnotte \cite[Chapter5]{Bonnotte} and gives the quantitative estimate
\begin{equation} \label{eq:sliced-comparison-bonnotte}
   \Wass_1(\mu,\nu) \lesssim R^{\frac{d}{d+1}}  \SlicedWass_1(\mu,\nu)^{\frac{1}{d+1}}.
\end{equation}
The second result concerns comparison with the stronger $\MaxSlicedWass_1$ sliced distance, for which Bobkov and G\"{o}tze in \cite{bobkov2024} proved 
\begin{equation} \label{eq:sliced-comparison-bobkov}
   \Wass_1(\mu,\nu) \lesssim R^{\frac{d}{d+2}}  \MaxSlicedWass_1(\mu,\nu)^{\frac{2}{d+2}},
\end{equation}
so that replacing $\SlicedWass_1$ by  $\MaxSlicedWass_1$ nearly doubles the exponent of Bonnotte.  The comparison \eqref{eq:sliced-comparison-bobkov} improves the exponents (also of the order of $\frac{2}{d}$ for large $d$) found in the 1980's by Hahn and Quinto in \cite{Hahn-Quinto}. The proof of \cite{bobkov2024} is simpler than the analysis of \cite{Hahn-Quinto} and also flexible enough to cover the case of measures with finite $q$-moments with $q>1$. Whether the exponent $\frac{2}{d+2}$ is sharp in \eqref{eq:sliced-comparison-bobkov} is not known but it is at least asymptotically sharp since \cite{bobkov2024} also explains why  in \eqref{eq:sliced-comparison-bobkov} the exponent cannot be better than $\frac{2}{d}$. 
Concerning the comparison \eqref{eq:sliced-comparison-bonnotte} between $\Wass_1$ and $\SlicedWass_1$, the
exponent $\frac{1}{d+1}$ of Bonnotte was thought to be optimal, see e.g. the remark following Theorem~3 in \cite{nadjahi2020statistical}. We prove in this paper that this is not true, and that the estimate \eqref{eq:sliced-comparison-bonnotte} can be slightly improved  to $\frac{1}{d}$ in odd dimensions and to $\frac{1}{d}$
up to a logarithmic correction in even dimensions. 
We also proved that (perhaps more surprisingly) the exponent $\frac{1}{d}$ is the sharp one. Let us finally mention the (widely open, in our opinion) case $p>1$. One can of course, again for compactly supported $\mu$ and $\nu$, directly derive quantitative comparisons between $\Wass_p$ by  $\SlicedWass_p$ from \eqref{eq:sliced-comparison-bonnotte} as Bonnotte did in \cite[Chapter5]{Bonnotte} obtaining an exponent $\frac{1}{p(d+1)}$  which is most likely far from optimal. In the present work, we will focus on the case $p=1$.
\hide{\edz{YW: I have not made any change to the narrative in the introduction. It needs to be done after the proofread of even dimensional case.}
}

\subsection*{Contributions} As already mentioned, our first contribution, to which Section \ref{sec-sharpcomp} is devoted, is an improvement of \eqref{eq:sliced-comparison-bonnotte} with an exponent $\frac{1}{d}$.  Moreover,  in Section \ref{sec-countergeneral}, we construct in any dimension $d$,  probability measures $p_\eps$ and $q_\eps$ supported on $B_1$ such that $\Wass_1(p_\eps, q_\eps) \gtrsim \eps$ and $\SlicedWass_1(p_\eps, q_\eps) \lesssim \eps^d$, which shows that the exponent $\frac{1}{d}$ is sharp.  Secondly, we extend our analysis in Section \ref{sec-kplaneslice} to the case of general sliced distances $\SlicedWass_1^k$ (see \eqref{eq:k-sliced} for the precise definition)  which involve average of Wasserstein distances between projections on $k$-planes of $\Rsp^d$. When $k=1$, this recovers the standard sliced-Wasserstein distance, while for $k=d-1$ this distance is closely related to the X-Ray transform. The main result of Section \ref{sec-kplaneslice} is that using $\SlicedWass_1^k$ instead of $\SlicedWass_1$ improves the exponent from $\frac{1}{d}$ to  $\frac{1}{d-k+1}$. 

\section{Sharp comparison between sliced Wasserstein and Wasserstein distances}\label{sec-sharpcomp}

\subsection{Notations and preliminaries}\label{sec-prel}

Let us first introduce some notations and recall some well-known results from optimal transport and Fourier analysis. We denote by $\Prob(\Rsp^d)$ the set of (Borel) probability measures on $\Rsp^d$ and for $p\geq 1$, we denote by  $\Prob_p(\Rsp^d)$ the set of probability measures with finite $p$-moment
\[\Prob_p(\Rsp^d)=\left\{\mu \in \Prob(\Rsp^d) \; : \;  M_p(\mu):=\int_{\Rsp^d} \vert x \vert^p \dd \mu(x)<+\infty\right\}.\] 
For $\mu$ and $\nu$ in $\Prob_p(\Rsp^d)$, the $p$-Wasserstein distance  $\Wass_p(\mu, \nu)$ between $\mu$ and $\nu$ is defined by
\[\Wass_p^p(\mu, \nu):=\inf_{\Gamma\in \Gamma(\mu, \nu)} \int_{\Rsp^d\times \Rsp^d} \vert x-y \vert^p \dd \gamma(x,y)\]
where $\Gamma(\mu, \nu)$ is the set of transport plans between $\mu$ and $\nu$ i.e. the set of probability measures on $\Rsp^d\times \Rsp^d$ having $\mu$ and $\nu$ as marginals. In the sequel, we shall focus on the case $p=1$, in this case the Kantorovich duality formula enables one to express $\Wass_1(\mu, \nu)$ as 
\begin{equation}\label{kanto-duality}
\Wass_1(\mu, \nu):=\sup \left\{ \int_{\Rsp^d} \phi \dd (\mu-\nu), \;   \phi \mbox{ $1$-Lipschitz} \right\}.
\end{equation}
Note first that $\Wass_1$ can be extended by homogeneity to pairs of nonnegative measures with the same total mass i.e. setting $\Wass_1(M \mu, M \nu)=M \Wass_1(\mu, \nu)$ for every $(M, \mu, \nu)\in \Rsp_+\times \Prob_1(\Rsp^d)^2$ and that $\Wass_1(\mu, \nu)$ is translation invariant i.e. only depends on $\mu-\nu$ (which is not the case for $p>1$) so that 
\begin{equation}
\Wass_1(\mu, \nu)=\Wass_1((\mu-\nu)^+, (\mu-\nu)^-)
\end{equation}
where $(\mu-\nu)^+$ and $(\mu-\nu)^-$ denote the positive and negative part of the signed measure $\mu-\nu$. This extension of $\Wass_1$ to signed measures will be useful in our examples showing sharpness (see Sections \ref{sec-countergeneral} 
). The Kantorovich duality \eqref{kanto-duality} formula of course implies the Kantorovich-Rubinstein inequality which states that for every Lipschitz function $\phi$, denoting by $\Lip(\phi)$ its Lipschitz constant, one has
\[\int_{\Rsp^d} \phi \dd (\mu-\nu) \leq  \Lip(\phi) \Wass_1(\mu, \nu).\]
If $d=1$, $\Wass_1$ has an appealingly simple expression (see Proposition 2.17 in \cite{Santambook}), 
\[\Wass_1(\alpha, \beta)= \Vert F_\alpha-F_\beta\Vert_{L^1(\Rsp)} =\int_\Rsp \vert F_{\alpha}(t)-F_\beta(t)\vert \dd t , \; \forall (\alpha, \beta)\in \Prob_1(\Rsp)\times \Prob_1(\Rsp)\]
where $F_\alpha$ is the cumulative distribution function of $\alpha$ i.e. $F_\alpha(t)=\alpha((-\infty, t])$ for every $t\in \Rsp$. Therefore for $(\mu, \nu)\in \Prob_1(\Rsp^d)^2$ denoting $(\mu^\theta, \nu^\theta)=(P_{\theta\#}\mu, P_{\theta\#}\nu)$, we have the following simple expression for $\SlicedWass_1(\mu,\nu)$
\begin{align}
\SlicedWass_1(\mu,\nu) &:=\int_{\theta \in\Sph^{d-1}}
\Wass_1(\mu^\theta,\nu^\theta) \dd \Haus^{d-1}(\theta) \nonumber\\
&=\int_{\theta \in\Sph^{d-1}} \Vert F_{\mu^\theta}-F_{\nu^\theta} \Vert_{L^1(\Rsp)} \dd \Haus^{d-1}(\theta) \label{eq:W11dcdf}.
\end{align}

We denote by $B_R$ the closed ball of radius $R$ centered at $0$. For $\mu$ and $\nu$ probability measures supported on $B_R$, the supremum in \eqref{kanto-duality} can be performed over $1$-Lipschitz functions on $B_R$ and is achieved by a (so-called Kantorovich potential) $\phi$ which may be extended to a $1$-Lipschitz function over $\Rsp^d$ with compact support included in $B_{3R}$. Indeed,  assuming that $\max_{B_{R}} \phi=0$ (which is without loss of generality since Kantorovich potentials are defined up to an additive constant) so that $\min_{B_{R}} \phi\geq -2R$ and then extending $\phi$ by infimal convolution  i.e. by $y \in \Rsp^d  \mapsto \min_{x\in B_{R}} \{\phi(x)+ \vert x-y\vert\}$, we see that $\phi$ is nonnegative outside $B_{3R}$ so taking the negative part of $\phi$ yields the desired $1$-Lipschitz extension that vanishes outside $B_{3R}$.

\smallskip

An important step for comparing $\Wass_1$ to $\SlicedWass_1$ will consist in using a known explicit representation (see Theorem 3.8 in \cite{Helgason2011})  of a given function $\psi$ (typically some smooth approximation of a Kantorovich potential in our case) in the form 
\begin{equation}\label{sphericalrep}
\psi(x) = \int_{\Sph^{d-1}} \psi^\theta(\sca{x}{\theta}) \dd \Haus^{d-1}(\theta), \; x\in \Rsp^d.
\end{equation}
For the sake of completeness, we recall how to obtain such an explicit representation using the Fourier transform. Let $\psi \in\Sc(\Rsp^d)$ (where $\Sc(\Rsp^d)$ is the Schwartz space of smooth  functions with rapidly decaying derivatives), and adopt the following convention for its Fourier transform
\[\Four{\psi}(y)=\hat{\psi}(y)= \int_{\Rsp^d} e^{-2i \pi \sca{x}{ y}} \psi(x) \dd x.\] 
Using Fourier inversion and spherical coordinates yields that for every $x\in \Rsp^d$
\[\psi(x)= \int_{\Rsp^d} e^{2i \pi \sca{x}{ y}} \hat{\psi}(y) \dd y=   \int_{\Sph^{d-1}}  \int_{\Rsp^+}  e^{2i \pi r \sca{x}{\theta}} \hat{\psi}(r \theta)r^{d-1} \; \dd r  \; \dd \Haus^{d-1}(\theta)\]
so that \eqref{sphericalrep} holds for $\psi^\theta$ given by
\begin{equation}\label{psiteta0}
 \psi^{\theta}(t) := \int_{\Rsp_+} \hat{\psi}(r\theta) e^{2i\pi r t} r^{d-1} \dd r= \frac{1}{(2i\pi)^{d-1}} \partial^{d-1}_t  \int_{\Rsp_+} \hat{\psi}(r\theta) e^{2i\pi r t}  \dd r.
 \end{equation}
 The Radon transform of $\psi$, $(\theta, t)\in \Sph^{d-1}\times \Rsp\to R^\theta\psi(t)$ is defined by
\begin{equation}\label{eq:defRadonT}
R^\theta \psi(t):=\int_{\theta^\perp} \psi(x + t\theta) \dd\Haus^{d-1}(x), \; \forall t\in \Rsp.
\end{equation}
Note that $R^\theta \psi \in \Sc(\Rsp)$ for every $\theta$ (and $R^\theta \psi \in \Class_c^\infty(\Rsp)$ if $\psi \in \Class_c^\infty(\Rsp^d)$) and that we have the so-called Fourier slicing formula: 
\[\hat{\psi}(r\theta)= \widehat{R^\theta \psi}(r).\]
 If $d$ is odd, the prefactor  $\frac{1}{(2i\pi)^{d-1}}$ in the right-hand side of \eqref{psiteta0} is real, taking the real part of the integral, we get
\begin{align*}
\psi^{\theta}(t) & =\frac{(-1)^{\frac{d-1}{2}}}{(2\pi)^{d-1}} \partial^{d-1}_t  \frac{1}{2} \Big( \int_{\Rsp_+} (\hat{\psi}(r\theta) e^{2i\pi r t} + \hat{\psi}(-r\theta) e^{-2i\pi r t}) \dd r\Big)\\
&=\frac{(-1)^{\frac{d-1}{2}}}{2^d \pi^{d-1}}   \partial^{d-1}_t \int_{\Rsp} e^{2i\pi r t} \widehat{R^\theta \psi}(r) \dd r= \frac{(-1)^{\frac{d-1}{2}}}{2^d \pi^{d-1}}  \partial^{d-1}_t  R^\theta \psi(t)
\end{align*}
where we used the fact that $\hat{\psi}(-r\theta)$ is conjugate to $\hat{\psi}(r\theta)$ in the first line, and Fourier inversion and the Fourier slicing formula in the second line. In the case of an even $d$, the prefactor in the right-hand side of \eqref{psiteta0} is purely imaginary, 
\[ \begin{split}
\psi^{\theta}(t)&= \frac{(-1)^{\frac{d}{2}}}    { (2\pi)^{d-1}}   \partial^{d-1}_t  \Big(i \int_{\Rsp_+}  \hat{\psi}(r\theta) e^{2i\pi r t}  \dd r \Big)\\
&=  \frac{(-1)^{\frac{d}{2}}}    { (2\pi)^{d-1}}     \partial^{d-1}_t  \Big(\frac{1}{2} \int_{\Rsp}  \hat{\psi}(r\theta) (-i\mathrm{sgn}(r)) e^{2i\pi r t}  \dd r\Big)\\
 &=  \frac{(-1)^{\frac{d}{2}}}   { 2^d \pi^{d-1}}     \partial^{d-1}_t  \Big( \int_{\Rsp}  \widehat{ R^\theta \psi}(r) (-i\mathrm{sgn}(r)) e^{2i\pi r t} \dd r \Big)
 \end{split}\]
 where we took the real part of the integral in the second line. To use a Fourier inversion again, we have to take into acount the odd sign prefactor which naturally makes a Hilbert transform appear. Recall that for $f\in \Sc(\Rsp)$, the Hilbert transform of $f$ is defined (in the sense of Cauchy principal value) as a convolution with the singular kernel $\frac{1}{s}$ i.e. by
\[Hf(t):= \frac{1}{\pi} \mathrm{pv} \int_{\Rsp} \frac{f(t-s)}{s} \dd s= \frac{1}{\pi} \lim_{\eps \to 0^+}  \int_{\Rsp\setminus [-\eps, \eps]} \frac{f(t-s)}{s} \dd s, \; t\in \Rsp.\]
Alternatively (and more conveniently for our context), $H$ can be defined in the Fourier domain as the inverse Fourier transform of $(-i \mathrm{sgn}) \Four$ (see chapter 5 of  \cite{grafakos2008fourier}) in the sense that, for $f\in \Sc(\Rsp)$ we have
\begin{equation}\label{HtransforminFourier}
\int_{\Rsp} \hat{f}(r) (-i \mathrm{sgn}(r)) e^{2 i \pi r t } \dd r= Hf(t). 
\end{equation}
For even $d$, we thus have 
\[\psi^{\theta}(t)=\frac{(-1)^{\frac{d}{2}}} { 2^d \pi^{d-1} }  \partial^{d-1}_t   (H R^\theta \psi)(t), \; \forall t\in \Rsp.\] 
So, to summarize, we have the representation \eqref{sphericalrep} of $\psi$ with $\psi^\theta$ given by
\begin{equation}\label{psitetaanydim}
\psi^\theta =\begin{cases} 
 \frac{(-1)^{\frac{d-1}{2}}} { 2^d \pi^{d-1} }  \partial^{d-1}_t  (R^\theta \psi)  \mbox{ if $d$ is odd},\\
\frac{(-1)^{\frac{d}{2}}} { 2^d \pi^{d-1} } \partial^{d-1}_t  (H R^\theta \psi) \mbox{ if $d$ is even}. 
\end{cases}
\end{equation}

\subsection{Control of  $\Wass_1$ by $\SlicedWass_1$.}

Our first result is an  improvement of the exponent  $\frac{1}{d+1}$ of Bonnotte:

\begin{Thm} \label{th:ComparisonW_SW}
For any dimension $d\geq 2$, there exists a constant $c_d$ such that for all $R>0$ and all $ \mu,\nu$ in $\Prob(\Rsp^d)$ supported on $B_R$, one has
\[ \Wass_1(\mu,\nu) \leq   
c_d R^{\frac{d-1}{d}} \SlicedWass_1(\mu,\nu)^{\frac{1}{d}}.\]

\hide{\edz{YW: statement of Theorem 2.1 is merged into one, with the power $\frac 1 d$. The log correction and its remarks are not deleted, just commented out, so that it is easier to check where the changes have been made.}
}
\end{Thm}

As we shall see  in Section \ref{sec-countergeneral}, in any dimension, the exponent $\frac{1}{d}$ cannot be improved so that the previous result is sharp. 
Theorem \ref{th:ComparisonW_SW} is a particular case of the more general Theorem \ref{th:ComparisonW_SWk} covering projections onto $k$-planes which will be the object of Section \ref{sec-kplaneslice} where a detailed proof will be given for all dimensions $d$ and $k$ with $1\leq k \leq d-1$. We will, however, give a short proof of Theorem \ref{th:ComparisonW_SW} in Section \ref{sec-sketch} for the odd-dimensional and in Section \ref{sec-proofeven}  for the even-dimensional case where some extra arguments are required to handle the Hilbert transform. 
Before doing so, let us consider some consequences of the previous theorem.

\smallskip

\textbf{Controlling $\Wass_1$ by finitely many one-dimensional marginals.} Let us remark that for $\mu$ and $\nu\in \Prob_1(\Rsp^d)$ and $\theta$, $\theta'$ in $\Sph^{d-1}$,  we obviously have
\[\Wass_1 ({P_\theta}_\#\mu, {P_{\theta'}}_\#\mu)\leq M_1(\mu) \vert \theta-\theta'\vert, \;  \Wass_1 ({P_\theta}_\#\nu, {P_{\theta'}}_\#\nu)\leq M_1(\nu) \vert \theta-\theta'\vert\]
which, together with the triangle inequality implies that $\Wass_1 ({P_\theta}_\#\mu, {P_\theta}_\#\nu)$ is an $(M_1(\mu)+M_1(\nu))$-Lipschitz function of $\theta$. Thanks to the Kantorovich-Rubinstein inequality, we thus deduce that for every probability measure $m$ on $\Sph^{d-1}$, we have
\[\SlicedWass_1(\mu,\nu)  \leq \int_{\Sph^{d-1}} \Wass_1( {P_\theta}_\#\mu, {P_\theta}_\#\nu) \dd m(\theta)+ (M_1(\mu)+M_1(\nu)) \Wass_1(\sigma_{d-1}, m),\]
where we recall that $\sigma_{d-1}$ denotes the uniform probability measure on $\Sph^{d-1}$. In particular, the previous inequality applies to empirical measures
\[\widetilde{m}_{N}(\widetilde\theta)=\frac{1}{N} \sum_{i=1}^N \delta_{\widetilde{\theta}_i}, \; \widetilde \theta=({\widetilde\theta}_1, \ldots, {\widetilde\theta}_N)\in (\Sph^{d-1})^N.\]
Combined with Theorem \ref{th:ComparisonW_SW}, this enables a quantitative control of $\Wass_1(\mu, \nu)$ in terms of the \emph{empirical sliced distance}
\[\widetilde{\SlicedWass}_1(\widetilde\theta, \mu,\nu):=\frac{1}{N} \sum_{i=1}^N  \Wass_1(  {P_{\widetilde{\theta}_i}}_\#\mu,  {P_{\widetilde{\theta}_i}}_\#\nu)
\]
and the \emph{quantification error}
\[\widetilde{\Delta}_N(\widetilde{\theta}):=\Wass_1(\sigma_{d-1},\widetilde{m}_{N}(\widetilde \theta) ).\]
For instance, for odd $d$ and for $\mu$ and $\nu$ supported on $B_R$, we have
\[\Wass_1(\mu,\nu) \leq c_d R^{\frac{d-1}{d}} \Big( \widetilde{\SlicedWass}_1(\widetilde \theta, \mu,\nu) +(M_1(\mu)+M_1(\nu)) \widetilde{\Delta}_N(\widetilde{\theta})\Big)^{\frac{1}{d}}.\]





\textbf{Relaxing the compact support assumption.} The statement of Theorem \ref{th:ComparisonW_SW} is restricted to compactly supported measures which is of course restrictive and exclude important cases such as Gaussians, we therefore propose the following variant for measures with exponential moments:

\begin{Cor}
Let $d$ be odd and $c_d$ denote the constant from Theorem \ref{th:ComparisonW_SW}. Then, for every $\mu$ and $\nu$ in $\Prob(\Rsp^d)$  such that for some $\alpha>0$, 
\[M_\alpha(\mu,\nu):=\int_{\Rsp^d} \vert x\vert e^{\alpha \vert x \vert} \dd (\mu+\nu)(x)\] 
is finite, we have
\[\Wass_1(\mu,\nu) \leq \SlicedWass_1(\mu,\nu)+ 2^{\frac{1}{d}} c_d \SlicedWass_1(\mu,\nu)^{\frac{1}{d}} \Big[\frac{1}{\alpha} \log \Big( \frac{M_\alpha(\mu,\nu)}{\SlicedWass_1(\mu,\nu)}\Big) \Big]^{\frac{d-1}{d}}.\]

\end{Cor}

\begin{proof}
For $R\geq 0$ denote by $\mathrm{proj}_{B_R}:\Rsp^d \to B_R$ the orthogonal projection onto $B_R$, namely
$$
\mathrm{proj}_{B_R}(x)=\left\{
\begin{array}{ll}
x& \text{if }|x|\leq R,\\
R\frac{x}{|x|} &\text{if }|x|>R,
\end{array}
\right.
$$
and set $\mu_R={\mathrm{proj}_{B_R}}_\#\mu$, $\nu_R={\mathrm{proj}_{B_R}}_\#\nu$. By the triangle inequality we have
\[\Wass_1(\mu,\nu) \leq \Wass_1(\mu,\mu_R) + \Wass_1(\mu_R,\nu_R)+ \Wass_1(\nu_R,\nu).\]
Also,  
\begin{align*}
\Wass_1(\mu,\mu_R) &\leq \int_{\Rsp^d} \vert x-\mathrm{proj}_{B_R}(x)\vert  \dd\mu(x)   \leq \int_{\Rsp^d\setminus B_R} |x| \dd\mu(x) \\
& \leq e^{-\alpha R} \int_{\Rsp^d}  \vert x\vert e^{\alpha \vert x \vert} \dd \mu(x)
\end{align*}
and similarly for $\Wass_1(\nu_R,\nu)$. Thus, using Theorem \ref{th:ComparisonW_SW}, we obtain
\begin{equation}\label{eq:firstboundexpo}
\Wass_1(\mu,\nu) \leq M_\alpha(\mu,\nu) e^{-\alpha R} +c_d R^{\frac{d-1}{d}} \SlicedWass_1(\mu_R,\nu_R)^{\frac{1}{d}}.
\end{equation}
Next, using the inequality $\SlicedWass_1\leq \Wass_1$, we deduce 
\begin{align*}
\SlicedWass_1(\mu_R,\nu_R) & \leq \SlicedWass_1(\mu,\nu)+\Wass_1(\mu,\mu_R)+\Wass_1(\nu,\nu_R)\\
& \leq \SlicedWass_1(\mu,\nu)+ M_\alpha(\mu,\nu) e^{-\alpha R}.
\end{align*}
Observing that $\SlicedWass_1(\mu,\nu) \leq \Wass_1(\mu,\nu) \leq \Wass_1(\mu,\delta_0)+ \Wass_1(\nu,\delta_0) \leq M_\alpha(\mu,\nu)$ (with a strict inequality unless $\mu=\nu=\delta_0$), we choose 
\[R=\frac{1}{\alpha} \log \Big( \frac{M_\alpha(\mu,\nu)}{\SlicedWass_1(\mu,\nu)}\Big)\]
so that $\SlicedWass_1(\mu_R,\nu_R) \leq 2 \SlicedWass_1(\mu,\nu)$. This, together with \eqref{eq:firstboundexpo} and the fact that $M_\alpha(\mu,\nu) e^{-\alpha R}=\SlicedWass_1(\mu,\nu)$, yields the desired inequality.

\end{proof}

\subsection{Proof of Theorem \ref{th:ComparisonW_SW} for odd $d$}\label{sec-sketch}

Let $\chi$ be a smooth probability density supported on the unit ball of $\Rsp^d$, and for $\lambda>0$ let $\chi_\lambda:= \lambda^{-d} \chi\Big( \frac{\cdot}{\lambda} \Big)$ be standard mollifiers and $\mu_\lambda:=\mu * \chi_\lambda$, $\mu_\lambda:=\mu * \chi_\lambda$ be smooth approximations of $\mu$ and $\mu$. The triangle inequality together with the fact that $\chi_\lambda$ is a probability measure supported on $B_\lambda$ imply that
\begin{equation}\label{trianglelambda}
\Wass_1(\mu,\nu) \leq   2 \lambda+\Wass_1(\mu_\lambda,\nu_\lambda) .
\end{equation}

The parameter $\lambda$ will be chosen later but we a priori already impose $\lambda \leq R$ so that $\mu_\lambda$, $\nu_\lambda$ are supported on $B_{2R}$. Let now $\phi$ be a Kantorovich potential between $\mu_\lambda$ and $\nu_\lambda$ i.e. a $1$-Lipschitz function $\phi$ on $B_{2R}$ such that
\[\Wass_1(\mu_\lambda,\nu_\lambda)=\int_{\Rsp^d} \phi \dd (\mu_\lambda-\nu_\lambda).\]
We may further assume (see Section \ref{sec-prel}) that $\phi$ is extended to a $1$-Lipschitz on $\Rsp^d$ which vanishes outside $B_{6R}$.
Let us define $\psi:=\chi_{\lambda}*\phi$ (we omit to write the explicit dependence in $\lambda$ to avoid too heavy notations) by construction $\psi$ is smooth, supported on $B_{7R}$ and $1$-Lipschitz. In particular,  we can represent it using the construction of Section \ref{sec-prel}, see \eqref{sphericalrep}, \eqref{eq:defRadonT} and \eqref{psitetaanydim} (odd $d$ case) which yields 
\[ \psi(x) = \int_{\Sph^{d-1}} \psi^\theta(\sca{x}{\theta}) \dd \Haus^{d-1}(\theta)\]
  with  (denoting by $c_d$ a dimensional constant which may change from a line to another)
 \[\psi^\theta(t)=c_d \partial_t^{d-1} \int_{\theta^\perp} \psi(x+t\theta) \dd\Haus^{d-1}(x).\]
Using \eqref{trianglelambda} together with the construction of $\phi$ and $\psi$, we obtain:
\begin{align}
\Wass_1(\mu,\nu)& \leq  2 \lambda+ \int_{\Rsp^d} \phi \dd (\mu_\lambda-\nu_\lambda)= 2 \lambda+ \int_{\Rsp^d} \psi \dd (\mu-\nu)\nonumber\\
&= 2 \lambda + \int_{\Sph^{d-1}} \int_{\Rsp} \psi^\theta \dd ({P_\theta}_\#\mu-{P_\theta}_\#\nu) \; \dd \Haus^{d-1}(\theta)\nonumber \\
&  \leq 2 \lambda + \int_{\Sph^{d-1}} \Lip(\psi^\theta) \Wass_1({P_\theta}_\#\mu, {P_\theta}_\#\nu)  \; \dd \Haus^{d-1}(\theta) \nonumber \\
& \leq 2 \lambda + c_d \sup_{\theta\in \Sph^{d-1}}  \Lip(\psi^\theta) \SlicedWass_1(\mu, \nu) \label{trianglelambda2}
\end{align}
where we used the Kantorovich-Rubinstein inequality in the third line. It remains to bound the Lipschitz constant of $\psi^\theta$ (which is smooth and compactly supported on $[-7R, 7R]$). For $t\in [-7R, 7R]$, denoting by $\partial_\theta:= \nabla \cdot \theta$ the partial derivative in direction $\theta$ and by $\partial^{(d-1)}_\theta \chi$ the $(d-1)$-th partial derivative of $\chi$ in direction $\theta$ we have 
\begin{align*} 
(\psi^{\theta})'(t) &=c_d \partial_t^{d-1}\int_{\theta^{\perp}}  \partial_\theta \psi(x+t\theta)\, \dd\Haus^{d-1}(x)\\
&=c_d \partial_t^{d-1}\int_{\theta^{\perp}}  \int_{\Rsp^d}  \partial_\theta \phi (y) \chi_{\lambda} (x+t\theta-y) \dd y \, \dd\Haus^{d-1}(x)\\
&= \frac{c_d}{ \lambda^{d-1}} \int_{\theta^{\perp}}  \int_{\Rsp^d}  \partial_\theta \phi (y)  \; \partial^{(d-1)}_\theta \chi \Big(\frac{x+t\theta-y}{\lambda}\Big) \lambda^{-d} \dd y \, \dd\mathcal H^{d-1}(x)\\
&= \frac{c_d}{ \lambda^{d-1}} \int_{\theta^{\perp}\cap B_{14R}}  \int_{\Rsp^d}  \partial_\theta \phi (y) \; \partial^{(d-1)}_\theta \chi \Big(\frac{x+t\theta-y}{\lambda}\Big) \lambda^{-d} \dd y \, \dd\mathcal H^{d-1}(x)
\end{align*}
where in the last line, we used the fact that the integrand vanishes unless $\vert x+t\theta-y \vert \leq \lambda \leq R$ and $\vert y \vert \leq 6 R$ which together with $\vert  t \vert\leq 7R$ implies that $\vert x\vert \leq 14R$. Using the fact that $\vert \partial_\theta \phi \vert \leq 1$ since $\phi$ is $1$-Lipschitz, we have
\begin{align*} 
\vert (\psi^{\theta})'(t) \vert & \leq   \frac{c_d}{ \lambda^{d-1}} \int_{\theta^{\perp}\cap B_{14R}}  \int_{\Rsp^d}  \vert \; \partial^{(d-1)}_\theta \chi \Big(\frac{x+t\theta-y}{\lambda}\Big) \vert \lambda^{-d} \dd y \, \dd\mathcal H^{d-1}(x)\\
&=  \frac{c_d}{ \lambda^{d-1}} \Haus^{d-1}(B_{14R}) \Vert \partial^{(d-1)}_\theta \chi \Vert_{L^1(\Rsp^d)} =c_d \frac{R^{d-1}}{\lambda^{d-1}}
\end{align*}
Recalling \eqref{trianglelambda2},  we thus get the bound 
\[ \Wass_1(\mu,\nu)\leq  2 \lambda +c_d R^{d-1} \lambda^{-(d-1)}\SlicedWass_1(\mu,\nu)+ \lambda.\]
Taking $\lambda = 2^{-\frac{1}{d}} R^{\frac{d-1}{d}} \SlicedWass_1(\mu,\nu)^{\frac{1}{d}}$  (so that $ \lambda \leq R$ since $\SlicedWass_1(\mu,\nu)\leq 2R$), in the previous inequality  obtain the desired inequality for the odd dimensional case
\begin{equation*}
    \Wass_1(\mu,\nu)\leq c_d R^{\frac{d-1}{d}}  
\SlicedWass_1(\mu,\nu)^{\frac{1}{d}}. 
\end{equation*}

\subsection{Proof of Theorem \ref{th:ComparisonW_SW} for even \texorpdfstring{$d$}{d}.}\label{sec-proofeven} We proceed as in the odd-dimensional case, the  extra difficulty we have to take into account is the presence of a Hilbert transform (in the $t$ variable) in the definition of $\psi^\theta$:
\[\psi^\theta(t)=c_d H \Big(\partial_t^{d-1} \int_{\theta^\perp} \psi(x+t\theta) \dd\Haus^{d-1}(x)\Big).\]
Note that $\psi^\theta$ is supported on $[-7 R, 7R]$, let us now estimate its Lipschitz constant, we have to bound
\[(\psi^{\theta})'(t) = \frac{c_d}{ \lambda^{d-1}} H \Big( \int_{\theta^{\perp}\cap B_{14R}}  \int_{\Rsp^d}  \partial_\theta \phi (y) \; \partial^{(d-1)}_\theta \chi \Big(\frac{x+t\theta-y}{\lambda}\Big) \lambda^{-d} \dd y \, \dd\mathcal H^{d-1}(x) \Big).
\]
At this stage, it is convenient to write $y=y_{\theta^\perp}+y_\theta \theta$ with $(y_{\theta^\perp},y_\theta)\in \theta^\perp \times \Rsp$ and to observe that the integrand inside $H$ vanishes unless $\vert x \vert \leq 14 R$ and $\vert x-y_{\theta^\perp} \vert \leq \lambda$. For $z\in \theta^\perp$ and $s\in \Rsp$, set
\[p_z(s):= \partial^{(d-1)}_\theta \chi \Big(z+s\theta\Big)\]
so that for every $z$, $p_z$ is supported on $[-1, 1]$, Lipschitz (uniformly in $z$) and $\int_{\Rsp} p_z(s) \dd s=0$ since $\chi$ is compactly supported. Thanks to 
Lemma \ref{expansion} and Remark \ref{rem:4.5} below, there is a constant $C$ for which
\begin{equation}\label{boundHpz} \vert H p_z(s) \vert \leq \frac{C}{1 + s^2}, \; \forall (z,s)\in \theta^\perp \times \Rsp.\end{equation}
We can then express $(\psi^{\theta})'(t)$ as
\[ \frac{c_d}{ \lambda^{d-1}} \int_{\theta^{\perp}\cap B_{14R}} \int_{\Rsp} \int_{\theta^\perp \cap B(x, \lambda) } \partial_\theta \phi (y) H p_{\frac{x-y_{\theta^\perp}}{\lambda}} \Big(\frac{t-y_\theta}{\lambda} \Big) \frac {\dd \Haus^{d-1}(y_{\theta^\perp}) }{\lambda^{d-1}} \frac{\dd y_\theta}{\lambda} \dd \Haus^{d-1}(x)\]
which, together with \eqref{boundHpz} and the fact that $\phi$ is $1$-Lipschitz yields
\[\vert (\psi^{\theta})'(t) \vert \leq \frac{ C R^{d-1}}{\lambda^{d-1}} \int_\Rsp \frac{1}{1+ \Big(\frac{t-y_\theta}{\lambda}\Big)^2} \frac{ \dd y_\theta}{\lambda}=\frac{ C R^{d-1}}{\lambda^{d-1}} \int_\Rsp \frac{1}{1+ y^2} \dd y. \]
Having bounded $\Lip (\psi^\theta)$ by $CR^{d-1} \lambda^{-d+1}$ and recalling \eqref{trianglelambda2}, we conclude exactly as in the odd-dimensional case.

\begin{Lem}\label{expansion}
Let $f:\mathbb R\to \mathbb R$ be a Lipschitz compactly supported function satisfying  $\int_{\mathbb R}f(s)s^{j-1}ds=0$ for $j=1,\ldots,k$. Then
$$
Hf(t)= O\left(\frac{1}{1+|t|^{k+1}}\right).
$$
\end{Lem}
\begin{proof}
Let us assume $A>0$ is a constant such that ${\rm supp}(f)\subset [-A,A]$.
Note that, if $t \in [-2A,2A]$, then $Hf(t)$ is also bounded (by the Lipschitz regularity of $f$).

On the other hand, for $|t| \geq 2A$ and $s \in [-A,A]$, then
$\frac{1}{t-s}=
\frac{1}{t}+O\left(\frac{s}{t^2}\right).
$
Therefore, since $\int_{\mathbb R}f(s)s^{j-1}ds=0$ for $j=1,\ldots,k$, and $f$ is compactly supported,
\begin{align*}
Hf(t)
=&\mathrm{pv} \int_{\mathbb R} \frac{f(s)}{t-s} ds\\
=& \frac{1}{t}\int_{\mathbb R}f(s)ds+
\cdots +\frac{1}{t^k}\int_{\mathbb R}f(s) s^{k-1}ds+
O\left(\frac{1}{t^{k+1}}|\int_{\mathbb R}f(s) s^{k}ds|\right)\\
=&O\left(\frac{1}{t^{k+1}}\right).
\end{align*}

\end{proof}

\begin{rem}\label{rem:4.5} 
If $f(s)=\partial_s^k g(s) $ for some smooth and compactly supported function $g(s)$, then $\int_{\mathbb R}f(s)s^{j-1}ds=0$ for all $j=1, \cdots, k$ by integration by parts.
\end{rem}

\subsection{The exponent $\frac{1}{d}$ is sharp}\label{sec-countergeneral} The aim of this section is to exhibit an example showing that for every $d\geq 3$, the exponent $\frac{1}{d}$ is sharp in the statement of Theorem \ref{th:ComparisonW_SW}. To start our construction, we consider  a radial probability measure $\mu$ concentrated on the unit ball of the hyperplane orthogonal to the vertical direction, $H = e_d^{\perp} \subseteq \Rsp^d$, and with density proportional to $(1-\vert x\vert^2)_+^{d+1}$ with respect to the $(d-1)$-Hausdorff measure. We denote by $f$ the density of the projection of $\mu$ on any line orthogonal to $e_d$ and passing through the origin.
One can check that for $\abs{t}\leq 1$, up to a  multiplicative constant that may change from line to line,
\begin{align*}
f(t)&= c \int_{H \cap \{x_1=t\} \cap B_1}  (1 - \abs{x}^2)^{d+1} \dd\Haus^{d-2}(x)\\
&= c \int_0^{\sqrt{1-t^2}} (1-t^2-r^2)^{d+1} r^{d-3} \dd r.
\end{align*}

Thus, $0\leq f(t) \leq c (1-t^2)^{d+1}$ and $f$ belongs to $\Class^{d+1}_c(\Rsp)$. Next, we consider a function $g:[-1,1]\to\Rsp$ which satisfies $\int_{-1}^1 g(s)ds=1$ and 
\begin{equation}\label{eq:ortho-g}
    \forall k\in\{1,\hdots,d+1\}, \qquad \int_{-1}^1 s^k g(s) \dd s=0,
\end{equation} 
and we set $g_\eps(s) = \frac{1}{\varepsilon} g(s/\varepsilon)$. Finally, we define $\mu_\eps$ as an average, weighted by $g_\eps$, of ``vertical" translations of $\mu$ (i.e. along $e_d$). More precisely, defining $\tau_v$ the translation by $v$, i.e. $\tau_v(x)=x+v$, we set
\[ \mu_\eps := \int_{-1}^1 g(s) \tau_{\eps s e_d \#} \mu \dd s 
\]
We note at this point that $\mu_\eps$ is not a probability measure, but a signed measure with total mass one. We first prove a lower bound on the $1$-Wasserstein distance between $\mu$ and $\mu_\varepsilon$ (still well-defined since $\mu$ and $\mu_\varepsilon$ have the same total mass as recalled in the Section \ref{sec-prel}). Using Kantorovich duality formula \eqref{kanto-duality}, we know that for any $1$-Lipschitz function $\phi$, $\Wass_1(\mu_\eps,\mu)\geq \int \phi \dd(\mu_\eps-\mu)$. Taking $\phi(x) = \abs{\sca{x}{e_d}}$, we get
\begin{equation}\label{boundbelowwmueps} \Wass_1(\mu_\eps,\mu)\geq \int \phi \dd(\mu_\eps-\mu) = \int_{-1}^1 g(s) \abs{\sca{\eps s e_d}{e_d}} \dd s = \eps \int_{-1}^1 \abs{s} g(s) \dd s. \end{equation}
Here we have used $\mu$ is supported on the hyperplane where $x_d=0$, so $\int \phi \dd \mu=0$. Let us also notice that we can choose $g$ such that $\int_{-1}^1 \abs{s} g(s) \dd s>0$.  For proving an upper bound on the sliced Wasserstein distance, it is convenient to introduce a notation: we will decompose any direction $\theta$ in $\Sph^{d-1}$ as $\theta = \cos\Theta \omega+ \sin\Theta e_d$ with $\omega \in e^{\perp}_d$ and $\Theta\in \Rsp$. Since on the hyperplane $H$, the projection on $\Rsp\theta$ is given by $P_\theta(x) = \cos\Theta \cdot \sca{\omega}{x}$, we have 
$P_{\theta\#}\mu = (t\mapsto t\cos(\Theta))_{\#} P_{\omega\#}\mu$. Thus, the density $f_\Theta$ of $P_{\theta\#}\mu:=\mu^\theta$ is 
$f_\Theta= \frac{1}{\vert \cos\Theta\vert}f(\frac{\cdot}{  \cos\Theta})$ and 

\[\mu_{\eps}^\theta:=P_{\theta\#}\mu_\eps= \int_{-1}^1 g(s) f_{\Theta}(\cdot-s\eps \sin\Theta) \dd s. \]
To estimate $\SlicedWass_1(\mu, \mu_\eps)$ we split the integral on the sphere in  two zones. In the first zone, which is the set of unit vectors $\theta$ close to the vertical direction $e_d$ i.e. $\vert \cos \Theta \vert \leq \eps$, we use the  crude upper bound $\Wass_1(P_{\theta\#}\mu, P_{\theta\#}\mu_\eps) \leq \Wass_1(\mu, \mu_\eps)\leq C \eps$ to obtain
\begin{align*}
    \int_{\{\theta \in \Sph^{d-1}\; : \; \abs{\cos \Theta}\leq \eps\}} \Wass_1(P_{\theta\#}\mu, P_{\theta\#}\mu_\eps) \dd \theta 
    &\leq  C\eps \int_{\{\theta \in \Sph^{d-1}\; : \;  \abs{\cos \Theta}\leq \eps\}} \dd \theta \\
    &\leq C \eps\eps^{d-1}=C \eps^d.
\end{align*} 

We now estimate $\Wass_1(\mu^\theta, \mu^\theta_\eps)$ in the zone where  $\abs{\cos(\Theta)} \geq \eps$. We recall that the Wasserstein-$1$ distance between probability measures on the real line is equal to the $L^1$ distance between their cumulative distribution functions. Introducing the cumulative distribution function $F$ of $f$, note that 
\[F_{\Theta}=\begin{cases} F(\cdot/ \cos\Theta) \mbox{ if $\cos \Theta >0$}\\ 1-F(\cdot /  \cos\Theta) \mbox{ if $\cos \Theta <0$} \end{cases}\]
is the cumulative distribution function of $f_\Theta$. We can therefore rewrite the Wasserstein-$1$ distance between the projections of $\mu$ and $\mu_\eps$ as
\begin{equation}\label{eq:W1-cdf-example}
\Wass_1(\mu^\theta, \mu^\theta_\eps) = \int_{-\vert \cos \Theta\vert - \eps}^{\vert \cos{\Theta} \vert+\eps} \left\vert \int_{-1}^1 (F_\Theta(t+\eps s \sin \Theta)-F_{\Theta}(t) )g(s) \dd s \right\vert \dd t.
\end{equation}
Thanks to a Taylor expansion and noting  that 
\[F^{(k)}_\Theta =  \mathrm{sign}(\cos\Theta) F^{(k)}(\cdot/\cos(\Theta))\cos(\Theta)^{-k}, \;  \; k\in\{1,\hdots,d+1\},\]
we get
\begin{align*}
\abs{F_{\Theta}(t+\eps s \sin \Theta)-F_{\Theta}(t) - \sum_{k=1}^{d+1} \frac{F_{\Theta}^{(k)}(t)}{k!} (\eps  s \sin \Theta)^k}
&\leq \frac{\Lip(F_\Theta^{(d+1)})}{(d+1)!} (\eps \vert s \vert \abs{\sin\Theta})^{d+2}\\
&\leq \frac{C(\eps \vert s\vert)^{d+2}}{\abs{\cos(\Theta)}^{d+2}}.
\end{align*}
 We now observe  that thanks to the orthogonality relations \eqref{eq:ortho-g}, the terms of degree $k\in\{1, \hdots, d+1\}$ in $s$ have zero mean when integrated against $g$.  Using the previous inequality to upper bound \eqref{eq:W1-cdf-example} then gives us
\[\Wass_1(\mu^\theta, \mu^\theta_\eps) \leq 2C(\abs{\cos\Theta}+ \eps) \frac{\eps^{d+2}}{\abs{\cos\Theta}^{d+2}} \leq
\frac{C\eps^{d+2}}{\abs{\cos\Theta}^{d+1}},\]
where we used the fact that $\Theta$ satisfies $\vert \cos(\Theta)\vert  \geq  \eps$ to get the last inequality. We now integrate this inequality on the part of the sphere where  $ \vert \cos \Theta \vert  \geq \eps$. Denoting $z$  the vertical component of $\theta \in\Sph^{d-1}$, we have  $z^2=1-\cos^2 \Theta\in [0, 1-\eps^2]$ so that (with a positive constant $C$ possibly changing from line to line)
\begin{align*}
\int_{\abs{\cos \Theta}\geq \eps}  \Wass_1(\mu^\theta, \mu^\theta_\eps) \dd \Haus^{d-1} (\theta)
&\leq C \eps^{d+2} \int_{\abs{\cos \Theta}\geq \eps} \frac{1}{\abs{\cos \Theta}^{d+1}} \dd\Haus^{d-1}(\theta)  \\
& = C\eps^{d+2} \int_{z^2 \leq 1-\eps^2} \frac{1}{(1-z^2)^{\frac{d+1}{2}}} (1-z^2)^{\frac{d-2}{2}} \dd z \\
&=  C \eps^{d+2} \int_0^{\sqrt{1-\eps^2}} \frac{1}{(1-z^2)^{3/2}} \dd z \\
&= C\eps^{d+2} \cdot \frac{x}{\sqrt{1-x^2}}\big|^{\sqrt{1-\eps^2}}_0\\
&\leq C  \eps^{d+2}\cdot \frac{1}{\eps}= C\eps^{d+1}.
\end{align*}
Summing the contributions of the two zones and recalling \eqref{boundbelowwmueps}, we find that for some $C>1$
\begin{equation}\label{bothsidesmueps}
\SlicedWass_1(\mu, \mu_\eps) \leq C \eps^{d}, \;  \Wass_1(\mu, \mu_\eps) \geq \frac{\eps}{C}.
\end{equation}
Note that since $g$ changes sign $\mu-\mu_\eps$ is a signed measure with zero average which we may write as $M_\eps (p_\eps- q_\eps)$ where $p_\eps$ and $q_\eps$ are probability measures (proportional to $(\mu-\mu_\eps)^+$ and  $(\mu-\mu_\eps)^-$ respectively and the normalizing constant $M_\eps$ is the total mass of $(\mu-\mu_\eps)^+$. Since  $(\mu-\mu_\eps)^+ \leq \mu+  \int_{-1}^1 \vert g(s)\vert \tau_{\eps s e_d \#} \mu \dd s$ we have, on the one hand, $M_\eps \leq M:=1 + \int_{-1}^1 \vert g(s)\vert ds$. On the other hand, since $\mu_\eps(H)=0$, we have
\[M_\eps \geq (\mu-\mu_\eps)(H) \geq \mu(H)=1.\] 
Hence we directly deduce from \eqref{bothsidesmueps} that the probability measures $p_\eps$ and $q_\eps$ (both supported on $B_1$) satisfy
\[ \Wass_1(p_\eps, q_\eps)= \frac{1}{M_\eps} \Wass_1(\mu, \mu_\eps) \geq \frac{\eps}{C M} \mbox{ and } \SlicedWass_1(p_\eps, q_\eps)=\frac{ \SlicedWass_1(\mu, \mu_\eps)}{M_\eps } \leq C  \eps^{d}.\]

\section{Comparison when slicing is performed by projections onto $k$-planes}\label{sec-kplaneslice}

We now consider the case where slicing is performed over $k$-planes and not simply lines. Given an integer $k$ between $1$ and $d-1$ we denote by $G(d, k)$ the Grassmannian of $k$-dimensional subspaces of $\Rsp^d$. Given $\xi \in G(d,k)$, we denote by $\xi^{\perp}$ its orthogonal and by $P_\xi$ the orthogonal projection onto $\xi$. It is well-know that there exists a unique invariant by rotation Borel probability measure on  $G(d,k)$, which we denote by $\sigma$ and view as a uniform probability measure on the Grassmanian. The sliced Wasserstein distance $\SlicedWass_1^k(\mu,\nu)$ is a natural generalization of $\SlicedWass_1$, which we define as the average over $G(d, k)$ of the Wasserstein distance $ \Wass_1({P_{\xi}}_\#\mu,{P_{\xi}}_\# \nu)$ between the projections. In other words, for $\mu$ and $\nu$ probability measures on $\Rsp^d$ with finite first moments, we set 
\begin{equation}\label{eq:k-sliced} 
\SlicedWass^k_1(\mu,\nu) := \int_{G(d, k)}
\Wass_1({P_{\xi}}_\#\mu,{P_{\xi}}_\# \nu) \dd\sigma(\xi).
\end{equation}

Our main result is the following:

\begin{Thm} \label{th:ComparisonW_SWk} 
For any dimensions $d>k\geq 1$, there exists a constant $C_{d,k}$ such that for all $R>0$ and all all $ \mu,\nu$ in $\Prob(\Rsp^d)$ supported on $B_R$, one has
\[ \Wass_1(\mu,\nu) \leq  
C_{d,k} R^{\frac{d-k}{d-k+1}} \SlicedWass_1^k(\mu,\nu)^{\frac{1}{d-k+1}} \]
\end{Thm}


This is of course a generalization of Theorem \ref{th:ComparisonW_SW} and the strategy of proof will be  the same, we will start by recalling well-known results from integral geometry and then will prove Theorem \ref{th:ComparisonW_SWk} first in the case where $d-k$ is even and then in the odd case, which requires a few extra arguments. In what follows $c_{d,k}$  (respectively $C_{d,k}$ or sometimes just $C$) will denote a  real (respectively positive) constant depending on $d$ and $k$ that may change from a line to another. 

\subsection{Some results from integral geometry}

Given $f\in \Sc(\Rsp^d)$, let us first recall that thanks to a well-known formula  due to Fuglede \cite{Fuglede1958} (also see the review \cite{Keinert1989}) 
\begin{equation}
\int_{\Rsp^d} f = c_{d,k} \int_{G(d,k)} \int_\xi \vert y\vert^{d-k} f(y) \dd\Haus^{k}(y)  \dd \sigma(\xi) \mbox{ with }  c_{d,k}=\frac{\pi^{\frac{d-k}{2}} \Gamma(\frac{k}{2})}{\Gamma(\frac{d}{2})}.
\end{equation}
Combining Fuglede's formula with a Fourier inversion, we see that for every $x\in \Rsp^d$, we have
\[\begin{split}
f(x)&=\int_{\Rsp^d} \hat{f}(y)e^{2 i \pi \sca{x}{y}} \dd y=c_{d,k} \int_{G(d,k)} \int_{\xi} \vert y\vert^{d-k} \hat{f}(y)e^{2 i \pi \sca{x}{ y}} \dd \Haus^{k}(y) \dd \sigma(\xi) \\
&=c_{d,k} \int_{G(d,k)} \int_{\xi} \vert y\vert^{d-k} \hat{f}(y)e^{2 i \pi \sca{P_\xi( x)} {y}} \dd \Haus^{k}(y) \dd \sigma(\xi) 
\end{split}\]
 which yields the representation 
 \begin{equation}\label{repfkp}
f(x)= \int_{G(d,k)} f^{\xi}(P_\xi(x)) \dd \sigma(\xi) 
\end{equation}
where we define for every $\xi \in G(d,k)$ and $z\in \xi$ 
\begin{equation}\label{defdexi}
 f^{\xi}(z)=c_{d,k} \int_{\xi} \vert y\vert^{d-k} \hat{f}(y)e^{2 i \pi \sca{z}{y}} \dd \Haus^{k}(y)= c_{d,k}(- \Delta)^{\frac{d-k}{2}} \int_{\xi} \hat{f}(y) e^{2 i \pi \sca{z} {y}} \dd\Haus^k(y)
 \end{equation}
where $(-\Delta)^{\frac{d-k}{2}}$ is the  $\frac{d-k}{2}$ power of the ($k$-dimensional) operator $(-\Delta)$ with respect to the $z$ variable belonging to the $k$-plane $\xi$. To be more precise, identifying $\xi$ with $\Rsp^k$, for $g\in \Sc(\xi)=\Sc(\Rsp^k)$, $(-\Delta)^{\frac{d-k}{2}} g$ is defined in Fourier via
\begin{equation}\label{defdeltafrac}
(-\Delta)^{\frac{d-k}{2}} g= \Four^{-1} (\vert \cdot \vert^{d-k} \Four(g)).
\end{equation}
In particular, if $d-k$ is even, $(-\Delta)^{\frac{d-k}{2}}$ coincides with the usual local differential operator of order $d-k$ obtained by raising the Laplacian to the power ${\frac{d-k}{2}}$, up to a normalizing multiplicative constant.

\smallskip

 Now let us observe that since $z\in \xi$ setting
\begin{equation}\label{defdefxiperp}
f_{\xi^\perp}(z):= \int_{\xi^\perp} f(z+v) \dd \Haus^{d-k}(v)
 \end{equation}
we have
\[\begin{split}
 \int_{\xi} \hat{f}(y) e^{2 i \pi \sca{z}{y}} \dd\Haus^k(y)&=\int_\xi   e^{2 i \pi z \cdot y}  \int_{\xi\times \xi^{\perp}} f(u+v) e^{-2 i \pi \sca{u}{y}}   \dd \Haus^{k}(u) \dd \Haus^{d-k}(v) \dd\Haus^k(y)\\
 &=  \int_\xi   e^{2 i \pi \sca{z}{y}}  \int_{\xi} e^{-2 i \pi \sca{u}{y}} f_{\xi^\perp}(u) \dd \Haus^{k}(u)  \dd\Haus^k(y)\\
 &=  \int_\xi   e^{2 i \pi \sca{z}{y}}  \widehat{f_\xi^\perp}(y) \dd\Haus^k(y)= f_{\xi^\perp}(z)
 \end{split}\]
where we used Fubini's Theorem  in the second line and Fourier inversion in the last one (observe that $f_{\xi^\perp} \in \Sc(\xi)$). Together with \eqref{defdexi}, this enables us to rewrite the representation \eqref{repfkp} in the more concise form 
\begin{equation}\label{kplaneinv}
f^{\xi}=c_{d,k}  (-\Delta)^{\frac{d-k}{2}} f_{\xi^\perp}, \; f_{\xi^\perp}(z):= \int_{\xi^\perp} f(z+v) \dd \Haus^{d-k}(v), \; \forall z\in \xi.
\end{equation}
For more  on representation and inversion formulas  arising in integral geometry, we refer the reader to the textbook of Helgason \cite{Helgason2011}.

\subsection{Proof of Theorem \ref{th:ComparisonW_SWk} for even $d-k$}\label{evendk}

We start exactly as in the proof of Section \ref{sec-sketch}, approximating $\mu$ and $\nu$  by $\mu_\lambda$, $\nu_\lambda$ by standard mollification using compactly supported mollifiers, and we recall \eqref{trianglelambda} .
The parameter $\lambda$ will be chosen later but we a priori already impose $\lambda \leq R$ so that $\mu_\lambda$, $\nu_\lambda$ are supported on $B_{2R}$. Let now $\phi$ be a Kantorovich potential between $\mu_\lambda$ and $\nu_\lambda$ i.e. a $1$-Lipschitz function $\phi$ on $B_{2R}$ such that
\[\Wass_1(\mu_\lambda,\nu_\lambda)=\int_{\Rsp^d} \phi \dd (\mu_\lambda-\nu_\lambda)\]
and we again assume that $\phi$ is extended to a $1$-Lipschitz on $\Rsp^d$ which vanishes outside $B_{6R}$.
Let us define $f:=\chi_{\lambda}*\phi$, by construction, $f$ is smooth, supported on $B_{7R}$ and $1$-Lipschitz so that it can be written in the form \eqref{repfkp}, with $f^\xi$ defined by \eqref{kplaneinv}.  Our task is to bound
\begin{align} 
\Wass_1(\mu_\lambda,\nu_\lambda)&=\int_{\Rsp^d} \phi \dd (\mu_\lambda-\nu_\lambda)=\int_{\Rsp^d} f \dd (\mu-\nu)=\int_{B_R} f \dd (\mu-\nu) \nonumber\\
&= \int_{G(d,k)} \int_{B_R} f^{\xi}(P_\xi(x))  \dd (\mu-\nu)(x) \dd \sigma(\xi) \nonumber \\
&= \int_{G(d,k)} \int_{\xi \cap B_R} f^{\xi}  \dd ({P_\xi}_\#\mu-{P_\xi}_\#\nu) \dd \sigma(\xi). \label{wasswithfxi}
\end{align}
Thanks to the Kantorovich-Rubinstein inequality, we have
\begin{align}\
\Wass_1(\mu_\lambda,\nu_\lambda) & \leq \int_{G(d,k)} \Lip(f^\xi, B_R) \; \Wass_1({P_\xi}_\#\mu, {P_\xi}_\#\nu) \dd \sigma(\xi) \nonumber\\
& \leq \sup_{\xi \in G(d,k)} \Lip(f^\xi, B_R)  \SlicedWass_1^k(\mu,\nu). \label{kanrub} 
\end{align}
Our goal is now to estimate the Lipschitz constant $\Lip(f^\xi)$ for some $k$-plane $\xi\in G(d,k)$. Using the representation \eqref{kplaneinv}, this amounts to estimating partial derivatives of the function $f_{\xi^\perp}$. In order to do so, we endow $\xi$ with a system of $k$ orthogonal coordinates. For $j\in \{1, \cdots, k\}$ and $\alpha\in \N^{k}$, for $z\in \xi$, we have by definition of $f=\phi *\chi_\lambda$ and recalling formula \eqref{kplaneinv} for $f_{\xi^\perp}$, we have
\[\partial_j \partial^{\alpha} f_{\xi^\perp}(z)=  \int_{\xi^\perp} \int_{\Rsp^d} \partial_j \phi(w)  \lambda^{-d- \vert \alpha\vert} \partial^\alpha \chi \Big(\frac{z+v-w}{\lambda}\Big) \dd \Haus^{d}(w) \dd \Haus^{d-k}  (v).\]
Moreover, since $f$ is supported on $B_{7R}$, $f_{\xi^\perp}$ is supported on $\xi \cap B_{7R}$. Assume then that $z\in \xi$ with $\vert z\vert \leq  7R$ and now observe that the integrand in the previous formula vanishes unless $\vert w \vert \leq 6R$ and $\vert z+v-w\vert \leq \lambda \leq R$. Since $\vert z\vert \leq  7 R$, these conditions imply that $\vert v\vert \leq 14 R$, so that we may restrict the integral to $B_{14R}$: 
\[\partial_j \partial^{\alpha} f_{\xi^\perp}(z)=  \lambda^{- \vert \alpha\vert} \int_{\xi^\perp\cap B_{14 R}}\int_{\Rsp^d} \partial_j \phi(w) \lambda^{-d} \partial^\alpha \chi \Big(\frac{z+v-w}{\lambda}\Big) \dd \Haus^{d}(w) \dd \Haus^{d-k}  \dd\Haus(v).\]
Recalling that $\phi$ is $1$-Lipschitz, we thus get
\begin{align}
\vert \partial_j \partial^{\alpha} f_{\xi^\perp}(z) \vert & \leq \lambda^{- \vert \alpha\vert} \int_{\xi^\perp\cap B_{14R}}\int_{\Rsp^d}  \lambda^{-d} \vert  \partial^\alpha \chi \Big(\frac{z+v-w}{\lambda}\Big) \vert \dd \Haus^{d}(w) \dd \Haus^{d-k} (v)\nonumber \\
&=  \lambda^{- \vert \alpha\vert} \int_{\xi^\perp\cap B_{14 R}}\int_{\Rsp^d}  \vert  \partial^\alpha \chi (w) \vert \dd \Haus^{d}(w) \dd \Haus^{d-k} (v)\nonumber \\
&= \lambda^{- \vert \alpha\vert}  \Haus^{d-k} (\xi^\perp\cap B_{14R}) \Vert  \partial^\alpha \chi\Vert_{L^1(\Rsp^d)} \nonumber  
  \end{align}
  so that $\partial^{\alpha} f_{\xi^\perp}$ is Lipschitz continuous on $B_{7R}$ (hence on $\xi$) with a Lipschitz constant which can be controlled as 
    \begin{equation}\label{lipalphaderiv} 
  \Lip(\partial^{\alpha} f_{\xi^\perp}, \xi\cap B_{7 R})= \Lip(\partial^{\alpha} f_{\xi^\perp}, \xi) \leq C R^{d-k} \lambda^{- \vert \alpha\vert}
  \end{equation}
  with $C$ a positive constant depending on $\vert \alpha\vert$, $d$ and $k$. 
    In particular, since $d-k$ is even, $f^\xi=c_{d,k} (- \Delta)^{\frac{d-k}{2}} f_{\xi^\perp}$ is $C(d, k) R^{d-k}  \lambda^{- d+k}$-Lipschitz on $B_R$:
 \[\Lip(f^\xi, B_R)  \leq C R^{d-k} \lambda^{-d+k}, \; \forall \xi\in G(d,k).\]
  Together with \eqref{kanrub} and \eqref{trianglelambda}, this yields
 \[\Wass_1(\mu_\lambda,\nu_\lambda) \leq C R^{d-k} \lambda^{-d+k}  \SlicedWass_1^k(\mu,\nu) + 2 \lambda\]
  taking $\lambda=2^{-\frac{1}{d-k+1}} R^{\frac{d-k}{d-k+1}}    \SlicedWass_1^k(\mu,\nu)^{\frac{1}{d-k+1}}$ in the previous inequality\footnote{Since $\SlicedWass_1^k(\mu,\nu) \leq 2R$, the choice of the prefactor $2^{-\frac{1}{d-k+1}}$ ensures that the condition $\lambda \leq R$ that we imposed in the beginning of the proof is satisfied.} exactly proves  Theorem \ref{th:ComparisonW_SWk} when $d-k$ is even.

\subsection{Proof of Theorem \ref{th:ComparisonW_SWk} for odd $d-k$}

If $d-k$ is odd, we have to handle the fact that in \eqref{kplaneinv} $(-\Delta)^{\frac{d-k}{2}}$ is  a nonlocal operator. We need some preliminary considerations about these operators and Riesz transforms.  For fixed $\xi\in G(d, k)$, let us identify $\xi$ with $\Rsp^k$ and for $\psi\in \Sc(\xi)=\Sc(\Rsp^k)$ and $j\in \{1, \ldots, k\}$, the $j$-th Riesz transform of $\psi$ is defined as the convolution of $\psi$ with the singular kernel $y_j\vert y \vert^{k+1}$ in the sense of Cauchy principal value, i.e. 
\begin{align}
 R_j \psi(z)&=\frac{\Gamma(\frac{k+1}{2})}{\pi^{\frac{k+1}{2}} } \mathrm{pv} \int_{\Rsp^k} \frac{(z_j-y_j)}{\vert z- y \vert^{k+1}} \psi(y) \dd \Haus^k(y)\nonumber\\
 &= \frac{\Gamma(\frac{k+1}{2})}{\pi^{\frac{k+1}{2}} }  \lim_{\eps \to 0^+}  \int_{\Rsp^k\setminus B_\eps} \frac{y_j}{\vert y \vert^{k+1}} \psi(z-y) \dd \Haus^k(y)      , \; z\in \Rsp^k. \label{defRiesz}
\end{align}
It is well-known (see Proposition 5.1.14 in \cite{grafakos2008fourier}) that $R_j$ is a generalization of the Hilbert Transform in several dimensions and that it can be characterized by its Fourier transform by 
\begin{equation}\label{rieszinfourier}
\Four (R_j \psi)(y)= -i \frac{y_j}{\vert y\vert} \hat{\psi}(y), \; y \in \Rsp^k, \; \psi \in \Sc(\Rsp^k). 
\end{equation}
Let $g\in \Sc(\Rsp^k)$, by \eqref{rieszinfourier}, we have
\[ \Four(R_j((-\Delta)^{\frac{d-k-1}{2}} \partial_j g))(y)=2 \pi \vert y\vert^{d-k-2} y_j^2 \; \Four(g)(y)\]
so that summing over $j$, we have
 \[ \Four(\sum_{j=1}^k R_j((-\Delta)^{\frac{d-k-1}{2}} \partial_j g))(y)=2 \pi \vert y\vert^{d-k} \Four(g)= 2 \pi \Four((-\Delta)^{\frac{d-k}{2}} g)\]
hence
\begin{equation}\label{decomplapl}
 (-\Delta)^{\frac{d-k}{2}} g=\frac{1}{2 \pi} \sum_{j=1}^k R_j (-\Delta)^{\frac{d-k-1}{2}} \partial_j g.
\end{equation}

We are now ready to prove the statement of Theorem \ref{th:ComparisonW_SWk} corresponding to the case where $d-k$ is odd.  Exactly as in Section \ref{evendk}, for $\lambda \in (0, R)$, we define $f=\phi * \chi_\lambda$ where $\phi$ a Kantorovich potential between $\mu_\lambda$ and $\nu_\lambda$, $1$-Lipschitz on $\Rsp^d$ and supported on $B_{6R}$. Thanks to \eqref{kplaneinv} and \eqref{decomplapl}, we can write $f$ as the average of $f^\xi \circ P_\xi$ over $G(d,k)$ with 
\[f^\xi:=c_{d,k}(- \Delta)^{\frac{d-k}{2}} f_{\xi^\perp}= \frac{c_{d,k}}{2 \pi} \sum_{j=1}^k   R_j (-\Delta)^{\frac{d-k-1}{2}} \partial_j f_{\xi^\perp}\]
that is 
\[f^\xi=\sum_{j=1}^k R_j \psi_j \mbox{ with } \psi_j:=\frac{c_{d, k}}{2\pi} (-\Delta)^{\frac{d-k-1}{2}} \partial_j f_{\xi^\perp}.\]
Our task is to bound the Lipschitz constant of $f^\xi$ (which is smooth and compactly supported on $\xi \cap B_{7R}$. Let $l\in \{1, \ldots, k\}$ our goal is to bound uniformly
\begin{equation}\label{derivfxi}
    \partial_l f^\xi=  \frac{c_{d,k}}{2 \pi} \sum_{j=1}^k   R_j (-\Delta)^{\frac{d-k-1}{2}} \partial^2_{jl} f_{\xi^\perp}.
    \end{equation}
First observe that for $z\in \xi$, we have 
\begin{align*}
   \partial_{l} f_{\xi^\perp} (z)=\int_{\xi^\perp} \int_{\Rsp^d} \partial_l \phi(w) \chi\Big(\frac{z+v-w}{\lambda} \Big ) \lambda^{-d} \dd \Haus^d(w) \dd \Haus^{d-k}(v)
\end{align*}
identifying $w\in \Rsp^d$ with its projections $(w_\xi, w_{\xi^\perp})\in \xi\times \xi^\perp$ onto $\xi$ and $\xi^\perp$ respectively,  this rewrites
\[ \partial_{l} f_{\xi^\perp} (z)=\int_{\xi^\perp} \int_{\Rsp^d} \partial_l \phi(w) \chi\Big(\frac{z-w_\xi}{\lambda}, \frac{v-w_{\xi^\perp}}{\lambda} \Big ) \lambda^{-d} \dd \Haus^d(w) \dd \Haus^{d-k}(v)\]
so that, if we set, for $x=(x_\xi, x_{\xi^\perp})\in \xi \times \xi^\perp)\in \xi \times \xi^\perp$
\[p_{j,x_\xi^\perp} (x_\xi):= (-\Delta)^{\frac{d-k-1}{2}} \partial_j \eta(x_\xi, x_\xi^\perp) \mbox{ (where } (-\Delta)^{\frac{d-k-1}{2}} \partial_j \mbox{ is with respect to  $x_\xi$)} \]
recalling  \eqref{derivfxi} and using the fact that $R_j$ commutes with dilations and translations, we have
\[\partial_l f^\xi (z)= \frac{C} {\lambda^{d-k}} \sum_{j=1}^k \int_{\xi^\perp} \int_{\Rsp^d} \partial_l \phi(w) R_j p_{j, \frac{v-w_{\xi^\perp}}{\lambda}} \Big(\frac{z-w_\eps}{\lambda} \Big ) \frac{ \dd \Haus^k(w_\xi)}{\lambda^k} \frac{ \dd \Haus^{d-k}(w_{\xi^\perp})}{\lambda^{d-k}} \dd \Haus^{d-k}(v)\]
since the integrand above vanishes unless $\vert v \vert \leq 14 R$ and $\vert w_\xi^\perp-v \vert \leq \lambda$, we obtain thanks to Lemma \ref{Rieszexpansion} and the fact that $\vert \partial_l \phi \vert \leq 1$
\begin{align*}\vert \partial_l f^\xi (z)\vert \leq \frac{C R^{d-k}}{\lambda^{d-k}} \sum_{j=1}^k \sup_{y \in \xi^\perp \cap B_1} \int_{\xi} \Big\vert R_j p_{j, y} \Big( \frac{z-w_\xi}{\lambda}\Big) \Big \vert\frac{\dd \Haus^{k}(w_\xi)}{\lambda^k}\\
\leq  \frac{C R^{d-k}}{\lambda^{d-k}} \int_{\xi} \frac{1}{1 + \vert w_\xi\vert^{k+1}} \dd \Haus^k(w_\xi).
\end{align*}
We thus have found a bound of the form $C R^{d-k}\lambda^{-d+k}$ for the Lipschitz constant of $f^\xi$ which enables us to conclude the proof exactly as in the  even $d-k$ case.
\begin{Lem}\label{Rieszexpansion}
Let $f:\mathbb R^k\to \mathbb R$ be a Lipschitz compactly supported function satisfying  $\int_{\mathbb R^k}f=0$, then for every $x\in \Rsp^k$ and $j=1, \ldots, k$
$$
\vert R_j f(x) \vert \leq \left(\frac{C}{1+|x|^{k+1}}\right).
$$
where $C$ is a constant depending on $k$, the Lipschitz constant of $f$ and its support.
\end{Lem}
\begin{proof}
Let $R>0$ be such that ${\rm supp}(f)\subset B_R$.
Note that, if $x \in B_{2R}$, then $R_j f(x)$ is can be directly bounded (by the Lipschitz regularity of $f$). On the other hand, for $|x| \geq 2R$ and $y \in B_R$, then it is easy to see that 
\begin{equation}\label{diffrkernel}
\Big \vert \frac{x_j-y_j}{\vert x-y \vert^{k+1}} - \frac{x_j}{\vert x\vert^{k+1}} \Big \vert \leq C \frac{\vert y \vert}{ \vert x \vert^{k+1}.}
\end{equation}
And since  $\int_{\Rsp^k}f=0$  and $f$ is  supported on $B_R$, we have
\begin{align*}
R_jf(x)
=&\mathrm{pv} \int_{\mathbb R^k} f(y)\frac{(x_j-y_j)}{\vert x-y\vert^{k+1}} \dd y\\
=& \mathrm{pv} \int_{\mathbb R^k} f(y) \Big( \frac{f(y)(x_j-y_j)}{\vert x-y\vert^{k+1}}-\frac{x_j}{\vert x \vert^{k+1}} \Big) \dd y
\end{align*}
together with \eqref{diffrkernel} this yields
\[\vert R_j f(x) \vert \leq \frac{C}{\vert x \vert^{k+1}} \int_{\Rsp^k} \vert y \vert \vert f(y) \vert \dd y \]
which proves the desired result.
\end{proof}

{\bf Acknowledgments:}  GC and AF  acknowledge the support of the Lagrange Mathematics and Computing Research Center. QM acknowledges the support of the Agence nationale
de la recherche, through the PEPR PDE-AI project (ANR-23-PEIA-0004). YW gratefully acknowledges that part of this work was carried out during a visit to Université Paris-Saclay.

\bibliographystyle{amsplain}   

\bibliography{references}

\providecommand{\bysame}{\leavevmode\hbox to3em{\hrulefill}\thinspace}
\providecommand{\MR}{\relax\ifhmode\unskip\space\fi MR }
\providecommand{\MRhref}[2]{%
  \href{http://www.ams.org/mathscinet-getitem?mr=#1}{#2}
}
\providecommand{\href}[2]{#2}
\begin{thebibliography}{10}

\bibitem{bobkov2024}
Sergey~G. Bobkov and Friedrich G\"{o}tze, \emph{Quantified cram\'er-wold continuity theorem for the kantorovich transport distance}, 2024.

\bibitem{Bonet2023}
Cl\'ement Bonet, Paul Berg, Nicolas Courty, Lucas~Drumetz François~Septier, and Minh-Tan Pham, \emph{Spherical sliced-{W}asserstein}, International Conference on Learning Representations, 2023.

\bibitem{Bonneel2015}
Nicolas Bonneel, Julien Rabin, Gabriel Peyr\'e, and Hanspeter Pfister, \emph{Sliced and {R}adon {W}asserstein barycenters of measures}, Journal of Mathematical Imaging and Vision \textbf{51} (2015), no.~1, 22--45.

\bibitem{Bonnotte}
Nicolas Bonnotte, \emph{Unidimensional and evolution methods for optimal transportation}, Phd. thesis (2013).

\bibitem{Cuturi2013}
Marco Cuturi, \emph{Sinkhorn distances: Lightspeed computation of optimal transport}, Advances in neural information processing systems, 2013.

\bibitem{deshpande2019max}
Ishan Deshpande, Yuan-Ting Hu, Ruoyu Sun, Ayis Pyrros, Nasir Siddiqui, Sanmi Koyejo, Zhizhen Zhao, David Forsyth, and Alexander~G Schwing, \emph{Max-sliced {W}asserstein distance and its use for {GAN}s}, Proceedings of the IEEE/CVF Conference on Computer Vision and Pattern Recognition, 2019, pp.~10648--10656.

\bibitem{fournier2015rate}
Nicolas Fournier and Arnaud Guillin, \emph{On the rate of convergence in {W}asserstein distance of the empirical measure}, Probability theory and related fields \textbf{162} (2015), no.~3-4, 707--738.

\bibitem{Fuglede1958}
Bent Fuglede, \emph{An integral formula}, Math. Scand. \textbf{6} (1958), 207--212. \MR{105724}

\bibitem{grafakos2008fourier}
Loukas Grafakos, \emph{Classical {F}ourier analysis}, second ed., Graduate Texts in Mathematics, vol. 249, Springer, New York, 2008. \MR{2445437}

\bibitem{Hahn-Quinto}
Marjorie~G. Hahn and Eric~Todd Quinto, \emph{Distances between measures from 1-dimensional projections as implied by continuity of the inverse {R}adon transform}, Z. Wahrscheinlichkeitstheorie verw. Gebiete \textbf{70} (1985), no.~3, 361--380.

\bibitem{Helgason2011}
Sigurdur Helgason, \emph{Integral geometry and {R}adon transforms}, Springer, New York, 2011. \MR{2743116}

\bibitem{Keinert1989}
Fritz Keinert, \emph{Inversion of {$k$}-plane transforms and applications in computer tomography}, SIAM Rev. \textbf{31} (1989), no.~2, 273--298. \MR{997459}

\bibitem{kitagawa2024}
Jun Kitagawa and Asuka Takatsu, \emph{Sliced optimal transport: is it a suitable replacement?}, 2024.

\bibitem{Kolouri2019a}
Soheil Kolouri, Kimia Nadjahi, Umut \c{S}im\c{s}ekli, Roland Badeau, and Gustavo~K. Rohde, \emph{Generalized sliced {W}asserstein distances}, The 33rd Conference on Neural Information Processing Systems, 2019.

\bibitem{Kolouri2019b}
Soheil Kolouri, Phillip~E. Pope, and Charles~E. Martin, \emph{Sliced-{W}asserstein auto-encoders}, In International Conference on Learning Representations, 2019.

\bibitem{Meunier2022}
Dimitri Meunier, Massimiliano Pontil, and Carlo Ciliberto, \emph{Distribution regression with sliced {W}asserstein kernels}, Proceedings of the 39th International Conference on Machine Learning, volume 162 of Proceedings of Machine Learning Research, vol. 162, 2022, p.~15501– 15523.

\bibitem{nadjahi2021sliced}
Kimia Nadjahi, \emph{Sliced-{W}asserstein distance for large-scale machine learning: theory, methodology and extensions}, Ph.D. thesis, Institut polytechnique de Paris, 2021.

\bibitem{nadjahi2020statistical}
Kimia Nadjahi, Alain Durmus, L{\'e}na{\"\i}c Chizat, Soheil Kolouri, Shahin Shahrampour, and Umut \c{S}im\c{s}ekli, \emph{Statistical and topological properties of sliced probability divergences}, Advances in Neural Information Processing Systems \textbf{33} (2020), 20802--20812.

\bibitem{park2023geometry}
Sangmin Park and Dejan Slep{\v{c}}ev, \emph{Geometry and analytic properties of the sliced wasserstein space}, arXiv preprint arXiv:2311.05134 (2023).

\bibitem{rabin2012wasserstein}
Julien Rabin, Gabriel Peyr{\'e}, Julie Delon, and Marc Bernot, \emph{Wasserstein barycenter and its application to texture mixing}, Scale Space and Variational Methods in Computer Vision: Third International Conference, SSVM 2011, Ein-Gedi, Israel, May 29--June 2, 2011, Revised Selected Papers 3, Springer, 2012, pp.~435--446.

\bibitem{Rustamov2020}
Subhabrata~Majumdar Raif~Rustamov, \emph{Intrinsic sliced {W}asserstein distances for comparing collections of probability distributions on manifolds and graphs}, arXiv:2010.15285, 2020.

\bibitem{Santambook}
Filippo Santambrogio, \emph{Optimal transport for applied mathematicians}, Progress in Nonlinear Differential Equations and their Applications, vol.~87, Birkh\"auser/Springer, Cham, 2015, Calculus of variations, PDEs, and modeling. \MR{3409718}

\end{thebibliography}


\end{document}